\documentclass[12pt]{amsart}
\usepackage{latexsym,amscd,amssymb}  
\theoremstyle{plain}
  \newtheorem{thm}{Theorem}[section]
  \newtheorem{prop}[thm]{Proposition}
  \newtheorem{lem}[thm]{Lemma}
  \newtheorem{cor}[thm]{Corollary}
\theoremstyle{definition}
  
  \newtheorem{exmp}[thm]{Example}
\theoremstyle{remark}
  \newtheorem{rem}[thm]{Remark}

\def\bR{{\mathbf R}}
\def\del{\operatorname{del}}
\def\11{\hat{{\bf 1}}}
\def\00{\hat{{\bf 0}}}

\def\rad{{\mathfrak r}}
\def\D{{\mathcal D}^\bullet}
\def\const{\underline{k}}

\def\wc{\omega^\bullet}

\def\II{I^\bullet}

\def\J{J^\bullet}
\def\M{M^\bullet}
\def\N{N^\bullet}
\def\P{P^\bullet}

\def\lk{\operatorname{lk}}
\def\str{\operatorname{st}}

\def\sp{\operatorname{Sp\acute{e}}}
\def\Sh{\operatorname{Sh}}
\def\Sq{\operatorname{Sq}(\Sigma)}

\def\Shc{\operatorname{Sh}_c}

\def\Qco{\operatorname{Qco}(Y)}

\def\TOR{\operatorname{Tor}_A}
\def\Or{{\mathcal {}or}}
\def\rH{\tilde{H}}
\def\rXi{\tilde{\chi}}
\def\cH{\mathcal{H}}
\def\m{{\mathfrak m}}
\def\NN{{\mathbb N}}
\def\const{\underline{k} }
\def\ZZ{{\mathbb Z}}
\def\RR{{\mathbb R}}

\def\mdZ{\operatorname{gr}_R}

\def\MdA{\operatorname{Gr}_A}

\def\mdZop{\operatorname{gr}_{R^!}}
\def\mod{\operatorname{mod}}
\def\vect{\operatorname{vect}}
\def\md{\mod_R}
\def\mde{\mod_\emptyset}

\def\mdop{\operatorname{mod}_{R^!}}

\def\Hom{\operatorname{Hom}}
\def\uHom{\underline{\Hom}}
\def\uExt{\underline{\Ext}}
\def\Dk{{\bf D}_k}
\def\Dko{{\bf D}_k^\op}
\def\Do{{\bf D}}

\def\op{{\sf op}}
\def\cHom{{\mathcal Hom} }

\def\Ext{\operatorname{Ext}}

\def\cmpl{{\sf c}}

\def\sp{\operatorname{Sp\acute{e}}}

\def\<{{\langle}}
\def\>{{\rangle}}

\def\tC{\tilde{C}}
\def\hSig{\hat{\Sigma}}
\def\h1{\hat{1}}

\def\too{\longrightarrow}

\def\Id{\operatorname{Id}}

\def\Proj{\operatorname{Proj}}

\def\cF{{\mathcal F}}
\def\cG{{\mathcal G}}
\def\ccF{{\mathcal F}^\bullet}
\def\ccG{{\mathcal G}^\bullet}

\numberwithin{equation}{section}
\begin{document}
\title[Incidence Algebras of Finite Regular Cell Complexes]
{Dualizing complex of the incidence algebra \\ 
of a finite regular cell complex}
\author{Kohji Yanagawa}
\address{Department of Mathematics, 
Graduate School of Science, Osaka University, Toyonaka, Osaka 
560-0043, Japan}
\email{yanagawa@math.sci.osaka-u.ac.jp}

\maketitle

\begin{abstract}
Let $\Sigma$ be a finite regular cell complex with $\emptyset \in \Sigma$, 
and regard it as a {\it poset} (i.e., partially ordered set) by inclusion.  
Let $R$ be the incidence algebra of the poset $\Sigma$ over a field $k$. 
Corresponding to the Verdier duality for constructible sheaves on $\Sigma$,   
we have a dualizing complex $\wc \in D^b(\mod_{R \otimes_k R})$
giving a duality functor from $D^b(\md)$ to itself.    
This duality is somewhat analogous to the Serre duality 
for projective schemes ($\emptyset \in \Sigma$ plays 
a similar role to that of ``irrelevant ideals").  
If $H^i(\wc) \ne 0$ for exactly one $i$, then the underlying 
topological space of $\Sigma$ is Cohen-Macaulay  
(in the sense of the Stanley-Reisner ring theory). 
The converse also holds if $\Sigma$ is a meet-semilattice as a poset
(e.g., $\Sigma$ is a simplicial complex). 
$R$ is always a Koszul ring with $R^! \cong R^\op$. 
The relation between the Koszul duality for $R$ and 
the Verdier duality is discussed. This result is a variant of a 
theorem of Vybornov. 
\end{abstract}

\section{Introduction}
Let $\Sigma$ be a finite regular cell complex, and 
$X := \bigcup_{\sigma \in \Sigma} \sigma$ its underlying topological space. 
The order given by $\sigma > \tau \stackrel{\text{def}}{\Longleftrightarrow} 
\bar{\sigma} \supset \tau$ makes $\Sigma$ a finite partially ordered set 
({\it poset}, for short). 
Here $\bar{\sigma}$ is the closure of $\sigma$ in $X$. 
Let $R$ be the incidence algebra of the poset $\Sigma$ over a field $k$. 
For a ring $A$, $\mod_A$ denotes the category of finitely generated 
left $A$-modules. 
In this paper, we study the bounded derived category $D^b(\md)$ 
using the theory of constructible sheaves (e.g., Verdier duality). 
For the sheaf theory, consult \cite{Iver, KS, Vyb}. 
We basically use the same notation as \cite{Iver}. 

Let $\Shc(X)$ be the category of $k$-constructible sheaves on $X$ 
with respect to the cell decomposition $\Sigma$. 
We have an exact functor $(-)^\dagger : \md \to \Shc(X)$.
For $M \in \md$, we have a natural decomposition 
$M = \bigoplus_{\sigma \in \Sigma} M_\sigma$ as a $k$-vector space. 
If $p \in \sigma \subset X$, the stalk $(M^\dagger)_p$ of $M^\dagger$ 
at the point $p$ is isomorphic to $M_\sigma$. 

Let $\Sigma' := \Sigma \setminus \emptyset$ be an induced subposet of 
$\Sigma$, and $T$ the incidence algebra of $\Sigma'$ over $k$. 
Then we have a category equivalence $\mod_T \cong \Shc(X)$, 
which is well-known to specialists (see for example 
\cite{Pls,Sh,Vyb}). But, in this paper, $\emptyset \in \Sigma$ 
plays a role.  Although $\md \not \cong \Shc(X)$, $\md$ has several 
interesting properties which $\mod_T$ does not possess. 
In some sense, $\emptyset$ is analogous to the ``irrelevant ideal" of 
a commutative noetherian homogeneous $k$-algebra. 

We have a left exact functor 
$\Gamma_\emptyset : \md \to \vect_k$ defined by 
$\Gamma_\emptyset(M) = \{ \, x \in M_\emptyset 
\mid Rx \subset M_\emptyset \, \}$. 
We denote its $i$th right derived functor by $H^i_\emptyset(-)$. 
For $M \in \md$, Theorem~\ref{cohomology} states that:  
$$H^i(X, M^\dagger) \cong H_\emptyset^{i+1}(M) \quad  
\text{for all $i \geq 1$,}$$  
$$0 \to H_\emptyset^0(M) \to M_\emptyset \to 
H^0(X, M^\dagger) \to H_\emptyset^1(M) \to 0 \quad \text{(exact)}.$$ 
Here $H^\bullet(X,M^\dagger)$ stands for the sheaf cohomology 
(c.f. \cite{KS, Iver}).

Let $A$ and $B$ be $k$-algebras. 
Recently, several authors study a dualizing complex  
$C^\bullet \in D^b(\mod_{A \otimes_k B})$ giving 
duality functors between $D^b(\mod_A)$ and $D^b(\mod_B)$. (Note that
if $M \in \mod_A$ and $N \in \mod_{A \otimes_k B}$ then 
$\Hom_A(M,N)$ has a left $B$-module structure.) 
In typical cases, it is assumed that $B = A^\op$. 
But, in this paper, from  Verdier's dualizing complex 
$\D_X \in D^b(\Shc(X))$ on $X$,  we construct a dualizing complex 
$\wc \in D^b(\mod_{R \otimes R})$ which gives the duality functor 
$\bR\Hom_R(-, \wc)$ from $D^b(\md)$ to itself. Theorem~\ref{Verdier} 
states that $$\bR\Hom_R(\M, \wc )^\dagger
\cong \bR\cHom((\M)^\dagger, \D_X)$$ in $D^b(\Shc(X))$ for all 
$\M \in D^b(\md)$. The dualizing complex $\wc$ satisfies the Auslander 
condition in the sense of \cite{YZ}. 

Corollary~\ref{Serre duality} states that 
$$\Ext^i_R(\M, \wc)_\emptyset \cong H^{-i+1}_\emptyset(\M)^\vee.$$ 
This corresponds to the (global) Verdier duality on $X$. 
But, since $H_\emptyset^i(-)$ can be seen as an analog of a local cohomology 
over a commutative noetherian homogeneous $k$-algebra, 
the above isomorphism can be seen as an imitation of the Serre duality. 
In Theorem~\ref{cohomology3} (1), $\emptyset \in \Sigma$ is also essential. 
It states that, for a simplicial complex $\Sigma$, $H^i(\wc) = 0$ 
for all $i \ne - \dim X$ if and only if $X$ is Cohen-Macaulay 
in the sense of the Stanley-Reisner ring theory. If we use the convention 
that $\emptyset \not \in \Sigma$, then the Cohen-Macaulay property 
can not be characterized in this way.

Under the assumption that a subset $\Psi$ of $\Sigma$ 
gives the open subset $U_\Psi := \bigcup_{\sigma \in \Psi} \sigma$ of $X$, 
Theorem~\ref{cohomology3} describes the cohomology 
$H^i(U_\Psi, M^\dagger|_{U_\Psi})$ using the duality functor 
$\bR \Hom_R(-, \wc)$. Note that the cohomology with compact support 
$H^i_c(U_\Psi, M^\dagger|_{U_\Psi})$ is much easier to treat 
in our context as shown in Lemma~\ref{open set}. 

We can regard $R$ as a graded ring in a natural way. 
Then $R$ is always Koszul, and the quadratic dual ring $R^!$ is isomorphic 
to the opposite ring $R^\op$ (Proposition~\ref{self dual}). 
Koszul duality (c.f. \cite{BGS}) gives an equivalence 
$D^b(\md) \cong D^b(\mod_{R^\op})$ of 
triangulated categories. The functors giving this equivalence 
coincide with the compositions of the duality 
functors $\bR\Hom_R(-, \wc)$ and $\Hom_k(-,k)$. 
This result is an ``augmented" version of Vybornov~\cite{Vyb}. 

It is well known that the  M\"obius function of a finite poset 
is a very important tool in combinatorics. In Proposition~\ref{Mobius}, 
generalizing \cite[Proposition~3.8.9]{StEn}, we describes the M\"obius 
function $\mu(\sigma, \h1)$ of the poset $\hSig := \Sigma \amalg \{ \h1 \}$ 
in terms of cohomology with compact support. 
As shown in \cite{Bj}, some finite posets arising from purely 
combinatorial/algebraic topics (e.g., Bruhat order)
are isomorphic to the posets of finite regular cell complexes. 
So the author expects that the results in the present paper will play a role 
in combinatorial study of these posets.

\section{Preparation}
A  {\it finite regular cell complex} (c.f. \cite[\S6.2]{BH} and 
\cite{CF}) is a non-empty topological space $X$ 
together with a finite set $\Sigma $ of subsets of $X$ such that the following 
conditions are satisfied: 
\begin{itemize}
\item[(i)] $\emptyset \in \Sigma $ and 
$X = \bigcup_{\sigma \in \Sigma } \sigma$; 
\item[(ii)] the subsets $\sigma \in \Sigma $ are pairwise disjoint;
\item[(iii)] for each $\sigma \in \Sigma $, $\sigma \ne \emptyset$, 
there exists a homeomorphism from an $i$-dimensional disc 
$B^i = \{ x \in \RR^i \mid ||x|| \leq 1 \}$ onto the closure 
$\bar{\sigma}$ of $\sigma$ which maps the open disc 
$U^i = \{ x \in \RR^i \mid ||x|| < 1 \}$ onto $\sigma$.  
\end{itemize}  
An element $\sigma \in \Sigma$ is called a {\it cell}. 
We regard $\Sigma$ as a poset by $\sigma > \tau 
\stackrel{\text{def}}{\Longleftrightarrow} \bar{\sigma} \supset \tau$.  
Combinatorics on posets of this type is discussed in \cite{Bj}. 
If $\sigma \in \Sigma $ is homeomorphic to $U^i$, 
we write $\dim \sigma = i$ and call $\sigma$ an {\it $i$-cell}. 
Here $\dim \emptyset = -1$. 
Set $d := \dim X = \max \{ \, \dim \sigma \mid \sigma \in  \Sigma  \, \}$. 

A finite simplicial complex is a primary example of finite regular cell 
complexes. When $\Sigma$ is a finite simplicial complex,   
we sometimes identify $\Sigma$ with the corresponding abstract simplicial 
complex. That is, we identify a cell $\sigma \in \Sigma$ with the set 
$\{ \, \tau \mid \text{$\tau$ is a 0-cell with $\tau \leq \sigma$} \, \}$. 
In this case, $\Sigma$ is a subset of the power set $2^V$, 
where $V$ is the set of the vertices (i.e., 0-cells) of $\Sigma$.  
Under this identification, for $\sigma \in \Sigma$, set 
$\str_\Sigma \sigma := \{ \, \tau \in \Sigma \mid 
\tau \cup \sigma \in \Sigma \, \}$ and 
$\lk_\Sigma \sigma := \{ \, \tau \in \str_\Sigma \sigma \mid 
\tau \cap \sigma = \emptyset \, \}$ to be subcomplexes of $\Sigma$.  

Let $\sigma, \sigma' \in \Sigma $. If $\dim \sigma =i+1$, 
$\dim \sigma' = i-1$ and $\sigma' < \sigma$, 
then there are exactly two cells $\sigma_1, \sigma_2 \in \Sigma $ between  
$\sigma'$ and $\sigma$. (Here $\dim \sigma_1 = \dim \sigma_2 = i$.) 
A remarkable property of a regular cell complex is the existence of an 
{\it incidence function} $\varepsilon$ (c.f. \cite[II. Definition~1.8]{CF}). 
The definition of an incidence function is the following. 
\begin{itemize}
\item[(i)] To each pair $(\sigma, \sigma')$ of cells, $\varepsilon$ assigns a 
number $\varepsilon(\sigma,\sigma') \in \{0, \pm 1\}$. 
\item[(ii)] $\varepsilon(\sigma,\sigma') \ne 0$ if and only if
$\dim \sigma' = \dim \sigma -1$ and $\sigma' < \sigma$.  
\item[(iii)] If $\dim \sigma = 0$, then $\varepsilon(\sigma,\emptyset) =1$. 
\item[(iv)] If $\dim \sigma =i+1$, $\dim \sigma' = i-1$ and 
$\sigma' < \sigma_1, \, \sigma_2 < \sigma$, $\sigma_1 \ne \sigma_2$, 
then we have 
$\varepsilon(\sigma, \sigma_1) \, \varepsilon(\sigma_1, \sigma') + 
\varepsilon(\sigma, \sigma_2) \, \varepsilon(\sigma_2, \sigma') =0$.
\end{itemize}
We can compute the (co)homology groups of $X$ using the cell decomposition 
$\Sigma$ and an incidence function $\varepsilon$.

\medskip

Let $P$ be a finite poset. The incidence algebra $R$ 
of $P$ over a field $k$ is the $k$-vector space with a basis 
$\{e_{x, \, y} \mid \text{$x, y \in P$ with $x \geq y$} \}$. 
The $k$-bilinear multiplication defined  by $e_{x, \, y} \, e_{z, \, w} = 
\delta_{y, \, z} \, e_{x, \, w}$ makes $R$ 
a finite dimensional associative $k$-algebra. Set $e_x := e_{x, \, x}$. 
Then $1 = \sum_{x \in P} e_x$ and $e_x \, e_y = \delta_{x,y} \, e_x$. 
We have $R \cong \bigoplus_{x \in P} R e_x$ as a left 
$R$-module, and each $R e_x$ is indecomposable. 

Denote the category of finitely generated left $R$-modules by $\md$. 
If $N \in \md$, we have $N = \bigoplus_{x \in P} N_x$ as a 
$k$-vector space, where $N_x := e_x N$.  Note that $e_{x, \, y} \, 
N_y \subset N_x$ and $e_{x,\, y} \, N_z  = 0$ for $y \ne z$. 
If $f: N \to N'$ is a morphism in $\md$, then $f(N_x) \subset N'_x$.  

For each $x \in P$, we can construct an indecomposable injective module 
$E_R(x) \in \md$. (When confusion does not occur, we simply 
denote it by $E(x)$.) Let $E(x)$ be the $k$-vector space with 
a basis $\{ e(x)_y \mid y \leq x \}$. 
Then we can regard $E(x)$ as a left $R$-module by 
\begin{equation}\label{envelope}
e_{z, \, w} \, e(x)_y  = 
\begin{cases}
e(x)_z & \text{if $y=w$ and $z \leq x$,} \\ 
0 & \text{otherwise.}
\end{cases}
\end{equation}
Note that $E(x)_y = k \, e(x)_y$ if $y \leq x$, and 
$E(x)_y = 0$ otherwise. An indecomposable injective in $\md$ is of the form 
$E(x)$ for some $x \in P$. Since $\dim_k R < \infty$, 
$\md$ has enough projectives and injectives. 
It is well-known that $R$ has finite global dimension. 

\medskip 
 
Let $\Sigma $ be a finite regular cell complex, and $X$ its underlying 
topological space.  We make $\Sigma $ a poset as above. 
In the rest of this paper, $R$ is the incidence algebra of 
$\Sigma $ over $k$.  For $M \in \md$, we have 
$M = \bigoplus_{\sigma \in \Sigma }M_\sigma$ 
as a $k$-vector space, where $M_\sigma := e_\sigma M$.

Let $\Sh(X)$ be the category of sheaves of finite dimensional 
$k$-vector spaces on $X$. We say $\cF \in \Sh(X)$ is a {\it constructible 
sheaf} with respect to the cell decomposition $\Sigma$, if 
$\cF|\sigma$ is a constant sheaf for all $\emptyset \ne \sigma \in \Sigma$. 
Here, $\cF|\sigma$ denotes the inverse image $j^* \cF$ of $\cF$ 
by the embedding map $j: \sigma \to X$.  Let $\Shc(X)$ be the full subcategory 
of $\Sh(X)$ consisting of constructible sheaves with respect to $\Sigma$.  
It is well-known that $D^b(\Shc(X)) \cong D^b_{\Shc(X)}(\Sh(X))$. 
(See \cite[Theorem~8.1.11]{KS}. There, it is assumed that $\Sigma$ 
is a simplicial complex. But this assumption is irrelevant. 
In fact, the key lemma \cite[Corollay~8.1.5]{KS} also holds for 
regular cell complexes. See also \cite[Lemma~5.2.1]{Sh}.) 
So we will freely identify these categories.

There is a functor $(-)^\dagger : \md \to \Shc(X)$ which 
is well-known to specialists (see for example \cite[Theorem~A]{Vyb}). 
But we give a precise construction here for the reader's 
convenience.   See \cite{Vyb, Y6} for detail. 

For $M \in \md$, set $$\sp(M) 
:= \bigcup_{\emptyset \ne \sigma \in \Sigma } \sigma \times M_\sigma.$$
Let $\pi : \sp(M) \to X$ be the projection map which sends $(p, m) \in 
\sigma \times M_\sigma \subset \sp(M)$ to $p \in \sigma \subset X$. 
For an open subset $U \subset X$ and a map $s: U \to \sp(M)$, 
we will consider the following conditions:  

\begin{itemize}
\item[$(*)$]  $\pi \circ s = \Id_{U}$ and $s_q = e_{\tau, \, \sigma} 
\cdot s_p$ for all $p \in \sigma$, $q \in \tau$ with $\tau \geq \sigma$. 
Here $s_p$ (resp. $s_q$) is the element of $M_\sigma$ 
(resp. $M_\tau$) with $s(p) = (p, s_p)$ (resp. $s(q) = (q, s_q)$).  
\item[$(**)$] There is an open covering $U = \bigcup_{\lambda \in \Lambda} 
U_\lambda$ such that the restriction of $s$ to $U_\lambda$ satisfies $(*)$ for 
all $\lambda \in \Lambda$. 
\end{itemize}

Now we define a sheaf $M^\dagger \in \Shc(X)$ from $M$ as follows. 
For an open set $U \subset X$, set 
$$M^\dagger(U):= 
\{ \, s \mid \text{$s: U \to \sp(M)$ is a map satisfying $(**)$} \,\}$$
and the restriction map $M^\dagger(U) \to M^\dagger(V)$ is the natural one. 
It is easy to see that $M^\dagger$ is a constructible sheaf. 
For $\sigma \in \Sigma$, let $U_\sigma := \bigcup_{\tau \geq \sigma} \tau$ 
be an open set of $X$. Then we have $M^\dagger(U_\sigma) \cong M_\sigma$. 
Moreover, if $\sigma \leq \tau$, then we have $U_\sigma \supset U_\tau$ and 
the restriction map $M^\dagger(U_\sigma) \to M^\dagger(U_\tau)$ 
corresponds to the multiplication map $M_\sigma \ni x \mapsto 
e_{\tau, \, \sigma} x \in M_\tau$. 
For a point $p \in \sigma$, the stalk $(M^\dagger)_p$ of $M^\dagger$ at 
$p$ is isomorphic to $M_\sigma$. 
This construction gives the functor $(-)^\dagger:\md \to \Shc(X)$. 
Let $0 \to M' \to M \to M'' \to 0$ be a complex in $\md$. 
The complex $0 \to (M')^\dagger \to M^\dagger \to (M'')^\dagger \to 0$ 
is exact if and only if 
$0 \to M'_\sigma \to M_\sigma \to M_\sigma'' \to 0$ is exact for all 
$\emptyset \ne  \sigma \in \Sigma$. 
Hence $(-)^\dagger$ is an exact functor. We also remark that 
$M_\emptyset$ is irrelevant to $M^\dagger$.  

For example, we have $E(\sigma)^\dagger \cong j_* \const_{\bar{\sigma}}$, 
where $j$ is the embedding map from the closure $\bar{\sigma}$ of $\sigma$ to 
$X$ and $\const_{\bar{\sigma}}$ is the constant sheaf on $\bar{\sigma}$. 
Similarly, we have $(R e_\sigma)^\dagger \cong h_! \const_{U_\sigma}$, 
where  $h$ is the embedding map from the open subset 
$U_\sigma = \bigcup_{\tau \geq \sigma} \tau$ to $X$. 

\begin{rem}
Let $\Sigma' := \Sigma \setminus \emptyset$ be an induced subposet of 
$\Sigma$, and $T$ its incidence algebra over $k$. 
Then we have a functor $\mod_T \to \Shc(X)$ defined by a similar way to 
$(-)^\dagger$, and it gives an equivalence $\mod_T \cong \Shc(X)$ 
(c.f. \cite[Theorem~A]{Vyb}). On the other hand, by virtue of 
$\emptyset \in \Sigma$, our $(-)^\dagger: \md \to \Shc(X)$ is neither 
full nor faithful. But we will see that $\md$ has several 
interesting properties which $\mod_T$ does not possess.
\end{rem}

For $M \in \md$, set 
$\Gamma_\emptyset(M) := \{ \, x \in M_\emptyset \mid 
Rx \subset M_\emptyset \, \}$. 
It is easy to see that $\Gamma_\emptyset(M)  \cong 
\Hom_R(k, M)$. Here we regard $k$ as a left $R$-module by 
$e_{\sigma, \, \tau} \, k = 0$ for all $e_{\sigma, \, \tau} \ne e_\emptyset$.
Clearly, $\Gamma_\emptyset$ gives a left exact functor 
from $\md$ to itself (or $\vect_k$). 
We denote the $i$th right derived functor of $\Gamma_\emptyset(-)$ 
by $H^i_\emptyset(-)$. In other words, $H_\emptyset^i(-) = 
\Ext_R^i(k, -)$. 

\begin{thm}[c.f. {\cite[Theorem~3.3]{Y6}}]\label{cohomology}
For $M \in \md$, we have an isomorphism 
$$H^i(X, M^\dagger) \cong H_\emptyset^{i+1}(M) \quad  
\text{for all $i \geq 1$},$$ and an exact sequence 
$$0 \to H_\emptyset^0(M) \to M_\emptyset \to 
H^0(X, M^\dagger) \to H_\emptyset^1(M) \to 0.$$
Here $H^\bullet(X,M^\dagger)$ stands for the cohomology 
with coefficients in the sheaf $M^\dagger$. 
\end{thm}

\begin{proof}
Let $\II$ be an injective resolution of $M$, 
and consider the exact sequence 
\begin{equation}\label{complexes}
0 \to \Gamma_\emptyset (\II) \to \II \to 
\II/ \Gamma_\emptyset(\II) \to 0
\end{equation}
of cochain complexes. Put $\J := \II/ \Gamma_\emptyset(\II)$. 
Each component of $\J$ is a direct sum of copies of 
$E(\sigma)$ for various $\emptyset \ne \sigma \in \Sigma $. 
Since $E(\sigma)^\dagger$ is the constant sheaf on $\bar{\sigma}$ 
which is homeomorphic to a closed disc, we have $H^i(X, E(\sigma)^\dagger) 
= H^i(\bar{\sigma}; k)= 0$ for all 
$i \geq 1$.  Hence $(\J)^\dagger \, (\, \cong (\II)^\dagger \,)$ 
gives a $\Gamma(X, -)$-acyclic resolution of $M^\dagger$. 
It is easy to see that $[\J]_\emptyset \cong \Gamma(X, (J^\bullet)^\dagger)$. 
So the assertions follow from \eqref{complexes}, since
$H^0(\II) \cong M$ and $H^i(\II) = 0$ for all $i \geq 1$. 
\end{proof}

\begin{rem}\label{dagger} 
(1) If $M_\emptyset = 0$, then we have 
$H^i(X, M^\dagger) \cong H_\emptyset^{i+1}(M)$ for all $i$.

(2) We regard a polynomial ring $S := k[x_0, \ldots, x_n]$ 
as a graded ring with $\deg(x_i) = 1$ for each $i$. 
Let $I \subset S$ be a graded ideal, and  set $A := S/I$. For a graded 
$A$-module $M$, we have the algebraic quasi-coherent sheaf $\tilde{M}$ 
on the projective scheme $Y := \Proj A$. It is well-known that 
$H^i(Y, \tilde{M}) \cong [H_\m^{i+1}(M)]_0$ 
for all $i \geq 1$,  and $$0 \to [H_\m^0(M)]_0 \to M_0 \to 
H^0(Y, \tilde{M}) \to [H_\m^1(M)]_0 \to 0 \qquad \text{(exact).}$$  
Here $H_\m^i(M)$ stands for the local cohomology module with support in 
the irrelevant ideal $\m := (x_0, \ldots, x_n)$, and 
$[H_\m^i(M)]_0$ is its degree 0 component ($H_\m^i(M)$ has a natural 
$\ZZ$-grading). See also Remark~\ref{section} (2) below.

(3) Assume that $\Sigma$ is a simplicial complex with $n$ vertices. 
The {\it Stanley-Reisner ring} $k[\Sigma]$ of $\Sigma$ is the quotient ring of 
the polynomial ring $k[x_1, \ldots, x_n]$ by the squarefree monomial ideal 
$I_\Sigma$ corresponding to $\Sigma$ (see \cite{BH,St} for details). 
In \cite{Y}, we defined {\it squarefree} $k[\Sigma]$-modules which are 
certain $\NN^n$-graded $k[\Sigma]$-modules. 
For example, $k[\Sigma]$ itself is squarefree. 
The category $\Sq$ of squarefree $k[\Sigma]$-modules is equivalent to $\md$ 
of the present paper (see \cite{Y5}). 
Let $\Phi: \md \to \Sq$ be the functor giving this 
equivalence. In \cite{Y6}, we defined a functor $(-)^+: \Sq \to \Shc(X)$. 
For example, $k[\Sigma]^+ \cong \const_X$. 
The functor $(-)^+$ is essentially same as the functor 
$(-)^\dagger: \md \to \Shc(X)$ of 
the present paper. More precisely, $(-)^\dagger \cong 
(-)^+ \circ \Phi$. For $M \in \md$, we have $H_\emptyset^i(M) \cong 
H_\m^i(\Phi(M))_0$. So the above theorem is a variation of 
\cite[Theorem~3.3]{Y6}. 
\end{rem}

\section{Dualizing complexes}
Let $D^b(\md)$ be the bounded derived category of $\md$. 
For $\M \in D^b(\md)$ and $i \in \ZZ$, 
$\M[i]$ denotes the $i$th translation of $\M$, 
that is, $\M[i]$ is the complex with $\M[i]^j = M^{i+j}$. 
So, if $M \in \md$, $M[i]$ is 
the cochain complex $\cdots \to 0 \to M \to 0 \to \cdots$, 
where $M$ sits in the  $(-i)$th position.

In this section, from Verdier's dualizing complex $\D_X \in D^b(\Shc(X))$, 
we construct a cochain complex $\wc$ of injective left 
$(R \otimes_k R)$-modules which gives a duality functor 
from $D^b(\md)$ to itself. 
Let $M$ be a left $(R \otimes_k R)$-module.  
When we regard $M$ as a left $R$-module via a ring homomorphism 
$R \ni x \mapsto x \otimes 1 \in R \otimes_k R$ 
(resp. $R \ni x \mapsto 1 \otimes x \in R \otimes_k R$ ), 
we denote it by ${}_R M$ (resp. $M_{R^\op}$). 

For $i \leq 1$, the $i$th component $\omega^i$ of $\wc$ has a $k$-basis 
$$\{ \, e(\sigma)^\tau_\rho \mid \sigma,\tau,\rho \in \Sigma , \, 
\dim \sigma  = -i, \, \sigma \geq \tau, \rho \, \},$$
and its module structure is defined by 
$$(e_{\sigma', \, \tau'} \otimes 1) \cdot e(\sigma)^\tau_\rho 
= \begin{cases} 
e(\sigma)^\tau_{\sigma'} & 
\text{if $\tau' = \rho$ and $\sigma' \leq \sigma$,} \\
0 & \text{otherwise,}
\end{cases}$$
and
$$
(1 \otimes e_{\sigma', \, \tau'}) \cdot e(\sigma)^\tau_\rho = 
\begin{cases}
e(\sigma)^{\sigma'}_\rho & 
\text{if $\tau'=\tau$ and $\sigma' \leq \sigma$,}\\
0 & \text{otherwise.}
\end{cases}$$

Then we have ${}_R (\omega^i) \cong (\omega^i)_{R^\op} \cong 
\bigoplus_{ \dim \sigma  = -i} E(\sigma)^{\mu(\sigma)}$ as left 
$R$-modules, where $\mu(\sigma) := \# \{ \tau \in \Sigma  \mid 
\tau \leq \sigma \}$. Note that $R \otimes_k R$ is isomorphic to the incidence 
algebra of the poset $\Sigma \times \Sigma $. For each $\sigma \in \Sigma $ 
with $\dim \sigma = -i$, set $I(\sigma)$ to be the subspace 
$\< \, e(\sigma)^\tau_\rho \mid \tau, \, \rho \leq \sigma \, \>$ 
of $\omega^i$. Then, as a left $R \otimes_k R$-module, 
$I(\sigma)$ is isomorphic to the 
injective module $E_{R \otimes_k R}( \, (\sigma,\sigma) \, )$, and 
$\omega^i \cong \bigoplus_{\dim \sigma  = -i} I(\sigma)$. 
Thus $\wc$ is of the form 
$$0 \to \omega^{-d} \to \omega^{-d+1} \to \cdots \to 
\omega^1 \to 0,$$
$$\omega^i = \bigoplus_{\substack{\sigma \in \Sigma \\ \dim \sigma = -i}}  
E_{R \otimes_k R} ( \, (\sigma, \sigma) \, ). $$
The differential of $\wc$ is given by 
$$\omega^i \ni e(\sigma)^\tau_\rho \longmapsto 
\sum_{\sigma' \geq \tau, \, \rho} 
\varepsilon(\sigma,\sigma') \cdot e(\sigma')^\tau_\rho \in \omega^{i+1}$$
makes $\wc$ a complex of $(R \otimes_k R)$-modules. 

Let $M \in \md$. Using the left $R$-module structure 
$I(\sigma)_{R^\op}$, $\Hom_R(M, {}_R I(\sigma))$ 
can be regarded as a left $R$-module again. 
Moreover, we have the following. 

\begin{lem}\label{injective}
For $M \in \md$, we have $\Hom_R(M,  {}_R I(\sigma)) \cong 
E(\sigma) \otimes_k (M_\sigma)^\vee$ as left $R$-modules. 
Here $(M_\sigma)^\vee$ is the dual vector space $\Hom_k(M_\sigma, k)$ of 
$M_\sigma$. 
\end{lem}

\begin{proof}
First, we show that if $M_\sigma = 0$ then 
$\Hom_R(M, {}_R I(\sigma)) =0$. Assume the contrary. If 
$0 \ne f \in \Hom_R(M, {}_R I(\sigma))$, there is some 
$x \in M_\tau$, $\tau < \sigma$, such that $f(x) \ne 0$. 
But we have $f(e_{\sigma, \, \tau} \, x) = e_{\sigma, \, \tau} \, f(x) \ne 0$. 
It contradicts the fact that $e_{\sigma, \, \tau} \, x \in M_\sigma = 0$. 

For a general $M \in \md$,  set $M_{\geq \sigma} = 
\bigoplus_{\tau \in \Sigma , \, \tau \geq \sigma} M_\tau$ 
to be a submodule of $M$. By the short exact sequence 
$0 \to M_{\geq \sigma} \to M \to M/M_{\geq \sigma} \to 0$, we have 
$$0 \to \Hom_R(M/M_{\geq \sigma}, {}_R I(\sigma)) \to 
\Hom_R(M, {}_R I(\sigma)) \to 
\Hom_R(M_{\geq \sigma}, {}_R I(\sigma)) \to 0.$$
Since $(M/M_{\geq \sigma})_\sigma = 0$, we have 
$\Hom_R(M, {}_R I(\sigma)) = 
\Hom_R(M_{\geq \sigma}, {}_R I(\sigma)).$ 
So we may assume that $M = M_{\geq \sigma}$. 
Let $\{ f_1,\ldots, f_n \}$ be a $k$-basis of $(M_\sigma)^\vee$. 
Since $({}_{R} I(\sigma))_\tau = 0$ for $\tau > \sigma$, 
$\Hom_R(M_{\geq \sigma}, {}_R I(\sigma))$ has a $k$-basis 
$\{ \, e(\sigma)_\sigma^\tau \otimes f_i 
 \mid \tau \leq \sigma, \, 1 \leq i \leq n \, \}$. 
By the module structure of $I(\sigma)_{R^\op}$, we have 
the expected isomorphism. 
\end{proof}

Since each ${}_R \omega^i$ is injective, 
$\Do(-) := \Hom^\bullet_R(-, {}_R\wc ) \cong \bR \Hom_R (-, {}_R\wc )$ 
gives a contravariant functor from $D^b(\md)$ to itself. 
In the sequel, we simply denote $\Hom_R(-, {}_R\omega^i )$ 
by $\Hom_R (-, \omega^i)$, etc. 

We can describe $\Do(\M)$ explicitly. 
Since $\omega^i \cong \bigoplus_{\dim \sigma  = -i} I(\sigma)$, we have 
$$\Hom_R(M, \omega^i) \cong 
\bigoplus_{\dim \sigma = -i} \Hom_R(M, I(\sigma)) 
\cong \bigoplus_{\dim \sigma  = -i} E(\sigma) \otimes_k  (M_\sigma)^\vee$$ 
for $M \in \md$ by Lemma~\ref{injective}. 
So we can easily check that $\Do(M)$ is of the form 
$$\Do(M) : 0 \too \Do^{-d}(M) \too \Do^{-d+1}(M) \too \cdots \too 
\Do^1(M) \too 0,$$ 
$$\Do^i(M) = \bigoplus_{\dim \sigma = -i} 
E(\sigma) \otimes_k (M_\sigma)^\vee.$$
Here the differential sends 
$e(\sigma)_\rho \otimes f \in E(\sigma) \otimes_k (M_\sigma)^\vee$ to  
$$\sum_{\tau \in \Sigma, \, \tau \geq \rho} 
\varepsilon(\sigma,\tau) \cdot
e(\tau)_\rho \otimes f(e_{\sigma, \tau} -) \ 
\in \bigoplus_{\dim \tau = \dim \sigma-1} 
E(\tau) \otimes_k (M_\tau)^\vee.$$
For a bounded cochain complex $\M$ of objects in $\md$, we have 
$$\Do^t(\M) = \bigoplus_{i-j = t} \Do^i(M^j) 
= \bigoplus_{-\dim \sigma -j = t}  
E(\sigma) \otimes_k (M^j_\sigma)^\vee,$$ 
and the differential is given by 
$$\Do^t(\M) \supset E(\sigma) \otimes_k (M^j_\sigma)^\vee  \ni 
x \otimes y \mapsto d(\, x \otimes y \, )+(-1)^t 
( x \otimes \partial^\vee(y)) \in \Do^{t+1}(\M), $$
where $\partial^\vee : (M^j_\sigma)^\vee \to (M^{j-1}_\sigma)^\vee$ 
is the $k$-dual of the differential $\partial$ of $\M$, 
and $d$ is the differential of $\Do(M^j)$. 

\medskip

Since the underlying space $X$ of $\Sigma $ is locally compact and 
finite dimensional, it admits Verdier's dualizing complex $\D_X \in 
D^b(\Sh(X))$ with the coefficients in $k$ (see \cite[V. \S2]{Iver}). 

\begin{thm}\label{Verdier}
For $\M \in D^b(\md)$, we have 
$$\Do(\M)^\dagger \cong \bR\cHom((\M)^\dagger, \D_X) \quad 
\text{in $D^b(\Shc(X))$}.$$ 
\end{thm}

\begin{proof}
An explicit description of $\bR\cHom((\M)^\dagger, \D_X)$ 
is given in the unpublished thesis \cite{Sh} of A. Shepard. 
When $\Sigma$ is a simplicial complex, this description is treated in 
\cite[\S 2.4]{Vyb}, and also follows from the author's previous paper 
\cite{Y6} (and \cite{Y5}). The general case can be reduced to 
the simplicial complex case using the barycentric subdivision. 

Shepard's description of $\bR\cHom((\M)^\dagger, \D_X)$ 
is the same thing as the above mentioned description of $\Do(\M)$
under the functor $(-)^\dagger$.  
\end{proof}

\begin{lem}\label{id for injectives}
For each $\sigma \in \Sigma$, the natural map 
$E(\sigma) \to \Do \circ \Do (E(\sigma))$ is an isomorphism in $D^b(\md)$. 
\end{lem}

\begin{proof}
We may assume that $\sigma \ne \emptyset$. Let 
$\Sigma|_\sigma := \{ \, \tau \in \Sigma \mid \tau \leq \sigma \, \}$
be a subcomplex of $\Sigma$. 
It is easy to see that $\Do(E(\sigma))_\emptyset$ is isomorphic 
to the chain complex $C_\bullet(\Sigma|_\sigma, k)$ of $\Sigma|_\sigma$. 
Thus $H^i(\Do(E(\sigma)))_\emptyset = \rH_{-i}( \bar{\sigma}; k)$ for all $i$, 
where $\rH_\bullet(\bar{\sigma};k)$ stands for the reduced homology group of 
the closure $\bar{\sigma}$ of $\sigma$. 
Hence $H^i(\Do(E(\sigma)))_\emptyset = 0$ for all $i$. 

By Theorem~\ref{Verdier} and the Poincar\'e-Verdier duality, we have 
$$\Do(E(\sigma))^\dagger \cong  \bR\cHom( j_*  \const_{\bar{\sigma}}, \D_X) 
\cong j_! \const_{\sigma}[\dim \sigma].$$  
Here $j: \sigma \to X$ is the embedding map.

Let $M$ be a simple $R$-module with $M = M_\sigma \cong k$. 
Combining the above observations, we have 
$\Do(E(\sigma)) \cong M[\dim \sigma]$. 
So $\Do \circ \Do(E(\sigma)) \cong \Do(M[\dim \sigma]) 
\cong E(\sigma)$, and the natural map 
$E(\sigma) \to \Do \circ \Do (E(\sigma))$ is an isomorphism. 
\end{proof}

\begin{thm}\label{Auslander}
(1) $\wc \in D^b(\mod_{R \otimes_k R})$ is a 
dualizing complex in the sense of \cite[Definition~1.1]{YZ}. 
Hence $\Do(-)$ is a duality functor from $D^b(\md)$ to itself.

(2) The dualizing complex $\wc$ satisfies the Auslander condition 
in the sense of \cite[Definition~2.1]{YZ}. That is, if we set 
$$j_\omega(M) := \inf \{ \, i \mid \Ext_R^i(M, \wc) \ne 0 \, \} 
\in \ZZ \cup \{ \infty \},$$ then, for all $i \in \ZZ$ and all $M \in \md$, 
any submodule $N$ of $\Ext_R^i(M,\wc)$ satisfies $j_\omega(N) \geq i$. 

\end{thm}

\begin{proof}
(1) The conditions (i) and (ii) of \cite[Definition~1.1]{YZ} obviously hold 
in our case. So it remains to prove the condition (iii). 
To see this, it suffices to show that the natural morphism 
$R \to \Do \circ \Do (R)$ is an isomorphism. 
But it follows from ``Lemma on Way-out Functors" 
(\cite[Proposition~7.1]{RD}) and Lemma~\ref{id for injectives}. 

(2) We may assume that $M \ne 0$. By the description of $\Do(M)$,  
we have $$j_\omega(M) = - 
\max\{ \, \dim \sigma \mid \sigma \in \Sigma, M_\sigma \ne 0 \, \}$$ 
and $\Ext_R^i(M, \wc)_\sigma = 0$ for $\sigma \in \Sigma$ with 
$\dim \sigma > -i$. Hence, any submodule  
$N \subset \Ext_R^i(M, \wc)$ satisfies $j_\omega(N) \geq i$. 
\end{proof}

\begin{cor}\label{Serre duality}
We have $\Ext^i_R(\M, \wc)_\emptyset \cong H^{-i+1}_\emptyset(\M)^\vee$ 
for all $i \in \ZZ$ and all $\M \in D^b(\md)$. 
\end{cor}

\begin{proof} 
Since $\Do \circ \Do(\M)$ is an injective resolution of $\M$, 
we have $\bR\Gamma_\emptyset(\M) = \Gamma_\emptyset(\Do \circ \Do(\M))$. 
By the structure of $\Do(-)$, we have  
$\Gamma_\emptyset(\Do \circ \Do(\M)) = (\Do(\M)_\emptyset)^\vee[-1]$. 
So we are done. 
\end{proof}

\section{Categorical Remarks}
For $M, N \in \md$ and $\sigma \in \Sigma$, 
set $\uHom_R(M,N)_\sigma := \Hom_R(M_{\geq \sigma}, N)$.  
We make $\uHom_R(M, N) := \bigoplus_{\sigma \in \Sigma} 
\uHom_R(M,N)_\sigma$ a left $R$-module as follows:  
For $f \in \uHom_R(M,N)_\sigma$ and a cell $\tau$ with $\tau \geq \sigma$, 
set $e_{\tau, \, \sigma} f$ to be the 
restriction of $f$ into the submodule $M_{\geq \tau}$ of $M_{\geq \sigma}$. 

\begin{lem}\label{injective2}
For $M \in \md$, we have $\uHom_R (M, E(\sigma)) \cong 
E(\sigma) \otimes_k (M_\sigma)^\vee$.
\end{lem}

\begin{proof}
Similar to Lemma~\ref{injective}. 
\end{proof}

If a complex $\M$ in $\md$ is exact, then 
so is $\uHom_R (\M, E(\sigma))$ by Lemma~\ref{injective2}. 
By the usual argument on double complexes, if  $\M$ is bounded and exact, 
and $\II$ is bounded and each $I^i$ is injective, 
then $\uHom^\bullet_R(\M, \II)$ is exact. 

\medskip

Note that $\Sigma$ is a {\it meet-semilattice} (see \cite[\S3.3]{StEn}) 
as a poset if and only if, for any two cells $\sigma, \tau \in \Sigma$ 
with $\bar{\sigma} \cap \bar{\tau}  \ne \emptyset$, there is a 
cell $\rho \in \Sigma$ with $\bar{\sigma} \cap \bar{\tau} = 
\bar{\rho}$. If $\Sigma$ is a 
simplicial complex, or more generally, a polyhedral complex, then 
it is a meet-semilattice.  If $\Sigma$ is a meet-semilattice, 
for two cells $\sigma, \tau \in \Sigma$, either 
there is no upper bound of $\sigma$ and $\tau$ (i.e., no cell 
$\rho \in \Sigma$ satisfies  $\rho \geq \sigma$ and $\rho \geq \tau$), 
or there is the least element $\sigma \vee \tau$ in 
$\{ \, \rho \in \Sigma \mid \rho \geq \sigma, \rho \geq \tau \, \}$
(c.f. \cite[Proposition~3.3.1]{StEn}).

Assume that $\Sigma$ is a meet-semilattice. 
Consider $\uHom_R(R e_\sigma, N)_\tau$ for $N \in \md$ and 
$\tau \in \Sigma$. If $\sigma \vee \tau$ exists, then 
we have $\uHom_R(R e_\sigma, N)_\tau = N_{\sigma \vee \tau}$. Otherwise, 
there is no upper bound of $\sigma$ and $\tau$, 
and $\uHom_R(R e_\sigma, N)_\tau = 0$. 
Hence the complex $\uHom_R(R e_\sigma, \N)$ is exact for an exact 
complex $\N$.  Hence if  $\N$ is bounded and exact, 
and $\P$ is bounded and each $P^i$ is projective, 
then $\uHom^\bullet_R(\P, \N)$ is exact. 

\medskip

By the above remarks, we have the following lemma 
(see \cite[I.1.10]{KS} for the derived functor of a bifunctor). 

\begin{lem}\label{RuHom} For $\M, \N \in D^b(\md)$, we have the following. 

(1) If $\II$ is an injective resolution of $\N$, then 
$$\bR \uHom_R(\M, \N) \cong \uHom_R^\bullet(\M, \II).$$

(2)  If $\Sigma$ is a meet-semilattice as a poset 
(e.g., $\Sigma$ is a simplicial complex), 
then $$\bR \uHom_R(\M, \N) \cong \uHom_R^\bullet(\P, \N)$$
for a projective resolution $\P$ of $\M$. 
\end{lem}

\begin{exmp}\label{2 disc}
The additional assumption in Lemma~\ref{RuHom} (2) is 
really necessary. That is, 
$\bR \uHom_R(\M, \N) \not \cong \uHom_R^\bullet(\P, \N)$ in general.   

For example, let $X$ be a closed 2 dimensional disc, 
and $\Sigma$ a regular cell decomposition of $X$ consisting 
of one 2-cell (say, $\sigma$), 
two 1-cells (say, $\tau_1, \tau_2$), and two 0-cells (say, $\rho_1, \rho_2$). 
Since $\rho_1 \vee \rho_2$ does not exist, 
$\Sigma$ is not a meet-semilattice.  

Let $N$ be a left $R$-module with $N = N_\sigma = k$. 
Then the injective resolution of 
$N$ is of the form $$\II : 0 \to E(\sigma) \to E(\tau_1) \oplus 
E(\tau_2) \to E(\rho_1) \oplus E(\rho_2) \to E(\emptyset) \to 0.$$ 
We have $\uHom_R(R e_{\rho_1}, E(\sigma))_{\rho_2} = 
\uHom_R(R e_{\rho_1}, E(\tau_1))_{\rho_2} 
= \uHom_R(R e_{\rho_1}, E(\tau_2))_{\rho_2} = k$ 
and $\uHom_R(R e_{\rho_1}, E(\rho_1))_{\rho_2} 
= \uHom_R(R e_{\rho_1}, E(\rho_2))_{\rho_2} = 0$. 
Thus $\uExt_R^1(R e_{\rho_1}, N)_{\rho_2} 
= H^1(\uHom(R e_{\rho_1}, \II))_{\rho_2} \ne 0$, 
while $R e_{\rho_1}$ is a projective module. 
\end{exmp}

\begin{prop}\label{uHom and Do} 
If $\M \in D^b(\Shc(X))$, then 
$\Do(\M) \cong \bR \uHom_R(\M, \Do(R e_\emptyset))$. 
\end{prop}

\begin{proof}
Since $\Do(R e_\emptyset)$ is of the form 
$0 \to D^{-d} \to D^{-d+1} \to \cdots \to D^1 \to 0$ with 
$D^i = \bigoplus_{\dim \sigma = -i} E(\sigma)$, the assertion follows from 
Lemmas~\ref{injective2} and \ref{RuHom}. 
\end{proof}

Since $(R e_\emptyset)^\dagger \cong \const_X$,  we have  
$\D_X \cong \Do(\const_X) \cong \Do(R e_\emptyset)^\dagger$ by 
Proposition~\ref{uHom and Do}. 

\medskip

If $\cF, \cG \in \Shc(X)$, then it is easy to see that 
$\cHom(\cF, \cG) \in \Shc(X)$. 
For $M, N \in \md$ and $\emptyset \ne \sigma \in \Sigma$, 
we have  
$\cHom (M^\dagger, N^\dagger) (U_\sigma) = \Hom_{\Sh(U_\sigma)} 
(M^\dagger|_{U_\sigma}, N^\dagger|_{U_\sigma}) 
= \Hom_R (M_{\geq \sigma}, N_{\geq \sigma})
= \Hom_R (M_{\geq \sigma}, N) = \uHom_R(M,N)_\sigma$. 
Hence $$\uHom_R(M,N)^\dagger \cong \cHom(M^\dagger, N^\dagger).$$

For $\ccF, \ccG \in D^b(\Shc(X))$, then it is known that 
$\bR \cHom (\ccF, \ccG) \in D^b(\Shc(X))$ 
(see \cite[Proposition~8.4.10]{KS}). Thus 
we can use an injective resolution of $\ccG$ in $D^b(\Shc(X))$ 
to compute $\bR \cHom (\ccF, \ccG)$. 
If $\II$ is an injective resolution of 
$\N \in D^b(\md)$, then $(\II)^\dagger$ is an injective resolution 
of $(\N)^\dagger$ in $D^b(\Shc(X))$. Hence we have the following.

\begin{prop}[{\cite[Theorem~5.2.5]{Sh}}]
If $\M, \N \in D^b(\md)$, then
$$\bR \uHom_R(\M, \N)^\dagger \cong \bR \cHom ((\M)^\dagger, (\N)^\dagger).$$
\end{prop}

By Lemma~\ref{RuHom} (2), if $\Sigma$ is a meeting-semilattice, 
then $\bR \cHom (\ccF, \ccG)$ for $\ccF, \ccG \in D^b(\Shc(X))$ 
can be computed 
using a projective resolution of $\ccF$ in $D^b(\Shc(X))$. 

\begin{rem}\label{section}
(1) Let $J$ be the left ideal of $R$ generated by 
$\{ \, e_{\sigma, \, \emptyset} \mid \sigma \ne \emptyset \, \}$. 
Note that $J^\dagger \cong \const_X$. 
Then we have that $\uHom_R(J, M)^\dagger \cong M^\dagger$ and 
$\uHom_R(J, M)_\emptyset \cong \Gamma(X, M^\dagger)$. 
Moreover, we have $\uExt_R^i(J, M) =\uExt_R^i(J, M)_\emptyset  
\cong H^i(X, M^\dagger)$ for all $i \geq 1$ 
by an argument similar to Theorem~\ref{cohomology}. 

(2) Let $\mde$ be the full subcategory 
of $\md$ consisting of modules $M$ with $M_\sigma = 0$ for all 
$\sigma \ne \emptyset$. 
Then $\mde$ is a {\it dense subcategory} of $\md$. That is, 
for a short exact sequence $0 \to M' \to M \to M'' \to 0$ in $\md$, 
$M$ is in $\mde$ if and only if $M'$ and $M''$ are in $\mde$. 
So we have the quotient category $\md/\mde$ by \cite[Theorem~4.3.3]{Po}. 
Let $\pi : \md \to \md/\mde$ be the canonical functor. It is easy to see that 
$\pi(M) \cong \pi(M')$ if and only if $M_{> \emptyset} \cong 
M'_{ > \emptyset}$. Moreover, we have $\Shc(X) \cong \md/\mde$. 

Let the notation be as in (1) of this remark. 
Then $\uHom_R(J,-)$ gives a functor 
$\eta: \md/\mde \to \md$ with $\pi \circ \eta = \Id$.  
Moreover, $\eta$ is a {\it section functor} 
(c.f. \cite[\S4.4]{Po}) and $\mde$ is a {\it localizing subcategory} of 
$\md$.

Let $A = \bigoplus_{i \geq 0}A_i$ be a commutative noetherian 
homogeneous $k$-algebra as in Remark~\ref{dagger} (2) and $\MdA$ the category 
of graded $A$-modules.  We say $M \in \MdA$ is a {\it torsion} 
module, if for all $x \in M$ there is some $i \in \NN$ with 
$A_{\geq i} \cdot x = 0$.  Let $\TOR$ be the full subcategory of $\MdA$ 
consisting of torsion modules.  Clearly, $\TOR$ is dense in $\MdA$. 
It is well-known that the category $\Qco$ of quasi-coherent sheaves 
on the projective scheme $Y := \Proj A$ is equivalent 
to the quotient category $\MdA/\TOR$, and  
we have the section functor $\Qco \to \MdA$ 
given by $\cF \mapsto \bigoplus_{i \in \ZZ} H^0(Y, \cF(i))$. 
So $\TOR$ is a localizing subcategory of $\MdA$. 
In this sense, our $\Shc(X) \cong \md/\mde$ is a small imitation 
of $\Qco \cong \MdA/\TOR$. 
\end{rem}

\section{Cohomologies of sheaves on open subsets}
Let $\Psi \subset \Sigma$ be an order filter of the poset $\Sigma$. That is, 
$\sigma \in \Psi$, $\tau \in \Sigma$, and $\tau \geq \sigma$ 
imply $\tau \in \Psi$. Then $U_\Psi := \bigcup_{\sigma \in \Psi} \sigma$ 
is an open subset of $X$. If $M \in \md$, $M_\Psi := 
\bigoplus_{\sigma \in \Psi} M_\sigma$ is a submodule of $M$. 
It is easy to see that $(M_\Psi)^\dagger \cong j_!j^*M^\dagger$, where 
$j : U_\Psi \to X$ is the embedding map. 
If $\Psi = \{ \tau \mid 
\tau \geq \sigma\}$ for some $\sigma \in \Sigma$, then 
$U_\Psi$ and $M_\Psi$ are denoted by 
$U_\sigma$ and $M_{\geq \sigma}$ respectively.

\begin{lem}\label{open set} 
Let $\Psi \subset \Sigma$ be an order filter with 
$\Psi \not \ni \emptyset$. 
Then we have the following isomorphisms for all $i \in \ZZ$ and $M \in \md$. 

(1) $H^{i+1}_\emptyset(M_\Psi) \cong H^i_c(U_\Psi, M^\dagger|_{U_\Psi})$ 
for all $i$. 

(2) $\Ext_R^i(M, \wc)_\sigma \cong 
H^{-i+1}_\emptyset(M_{\geq \sigma})^\vee \cong 
H^{-i}_c(U_\sigma, M^\dagger |_{U_\sigma})^\vee$ for all 
$\emptyset \ne \sigma \in \Sigma$. 
\end{lem}

\begin{proof} 
(1) We have $H^{i+1}_\emptyset(M_\Psi) \cong 
H^i(X, (M_\Psi)^\dagger) \cong H^i(X, j_! j^* M^\dagger) \cong 
H^i_c(U_\Psi, M^\dagger|_{U_\Psi})$.
Here, by Remark~\ref{dagger} (1),  the first isomorphism holds 
even if $i=0$. 

(2) By the description of $\Do(M)$, we have $\Do(M)_\sigma \cong 
\Do(M_{\geq \sigma})_\emptyset$. Hence we have 
$$\Ext_R^i(M, \wc)_\sigma \cong \Ext_R^i(M_{\geq \sigma}, \wc)_\emptyset 
\cong H^{-i+1}_\emptyset(M_{\geq \sigma})^\vee \cong 
H^{-i}_c(U_\sigma, M^\dagger|_{U_\sigma})^\vee.$$  
Here the second isomorphism follows from Corollary~\ref{Serre duality}. 
\end{proof}

\begin{prop}
For any $\sigma \in \Sigma$, 
$\Do(R e_\sigma)^\dagger \cong \bR j_* \D_{U_\sigma}$
where $j: U_\sigma \to X$ is the embedding map. In particular, 
$\Do(R e_\emptyset)^\dagger \cong \D_X$. 
\end{prop}

\begin{proof}
Set $U := U_\sigma$. Since $(R e_\sigma)^\dagger \cong j_! \const_U$, 
we have  
\begin{eqnarray*}
\Do(R e_\sigma)^\dagger 
&\cong& \bR \cHom (\, j_! \const_U, \, 
\D_X \, ) \quad \quad \text{(by Theorem~\ref{Verdier})}\\
&\cong& \bR j_* \bR\cHom (\,  \const_U, \, j^* \D_X \, ) 
\quad \quad \text{(by \cite[VII. Theorem~5.2]{Iver})}\\
&\cong& \bR j_* \bR\cHom(\, \const_U, \, \D_U \, )  \cong \bR j_* \D_U. 
\end{eqnarray*}
\end{proof}

Motivated by Lemma~\ref{open set},  
we will give a formula on the ordinal (not compact support) cohomology 
$H^i(U_\Psi, M^\dagger|_{U_\Psi})$. 

\begin{thm}\label{cohomology3}
Let $\Psi \subset \Sigma$ be an order filter with 
$\Psi \not \ni \emptyset$. We have 
$$H^i(U_\Psi, M^\dagger|_{U_\Psi}) \cong[ \, \Ext^i_R(\, \Do(M)_\Psi, 
\, \wc \, ) \, ]_\emptyset$$ for all $i \in \NN$ and $M \in \md$.
\end{thm}

\begin{proof}
For the simplicity, set $U := U_\Psi$. Let $\cF^\bullet \in D^b(\Sh(U))$.  
Taking a complex in the isomorphic class of $\cF^\bullet$, 
we may assume that each component
$\cF^i$ is a direct sum of sheaves of the form $h_! \const_V$, where 
$V$ is an open subset of $U$ with the embedding map $h: V \to U$ 
(see \cite[II. Proposition~2.4]{Iver}). 
Since each component ${\mathcal D}_U^i$ of $\D_U$ is an injective sheaf, 
$h^*{\mathcal D}_U^i$ is also injective by \cite[II. Corollary~6.10]{Iver},  
and we have $$\cHom (h_! \const_V, {\mathcal D}_U^i) \cong 
\bR h_* \bR \cHom (\const_V, h^*{\mathcal D}_U^i) \cong 
\bR h_*(h^*{\mathcal D}_U^i) \cong 
h_*h^*{\mathcal D}_U^i $$ by \cite[VII, Theorem~5.2]{Iver}.  
Since the sheaf $h_* h^*{\mathcal D}_U^i$ is flabby, 
$\cHom^\bullet(\cF^\bullet, \D_U)$ is a complex of flabby sheaves. 
Hence we have 
\begin{eqnarray*}
\Ext^i_{\Sh(U)}( \, \cF^\bullet, \D_U \, ) 
&\cong& H^i( \, \Gamma(U, \bR \cHom^\bullet
( \, \cF^\bullet, \, \D_U \, ) \, )\\ 
&\cong& \bR^i\Gamma( \, U, 
\bR \cHom( \, \cF^\bullet, \, \D_U \, ) \, ).
\end{eqnarray*}   

Since $\bR\cHom(\, \bR\cHom(M^\dagger|_U,\D_U), \, \D_U \, ) 
\cong M^\dagger|_U$ in $D^b(\Sh(U))$,  
we have 
\begin{eqnarray*}
H^i(U, M^\dagger|_U)
&\cong& \bR^i \Gamma(\, U, \, \bR\cHom
(\, \bR\cHom (M^\dagger|_U,\D_U), \, \D_U) \, ) \\  
&\cong& \Ext^i_{\Sh(U)}(\, \bR\cHom (M^\dagger|_U,\D_U), \, \D_U \, ) \\
&\cong& \bR^{-i} \Gamma_c(U, \bR\cHom (M^\dagger|_U, \D_U))^\vee 
\qquad \ (\text{by \cite[V, Theorem~2.1]{Iver}})\\
&\cong& \bR^{-i} \Gamma_c ( U, \bR\cHom (M^\dagger,\D_X)|_U)^\vee \\
&\cong& \bR^{-i} \Gamma_c ( U, \Do(M)^\dagger |_U )^\vee \\ 
&\cong& \bR^{-i+1} \Gamma_\emptyset (U, \Do(M)_\Psi)^\vee 
\qquad \ (\text{by Lemma~\ref{open set})} \\ 
&\cong& ( \, \Ext^i_R( \, \Do(M)_\Psi, \, \wc \, )_\emptyset \, )
\qquad \ (\text{by Corollary~\ref{Serre duality}}). 
\end{eqnarray*} 
\end{proof}

\begin{exmp}
Assume that $X$ is a $d$-dimensional manifold (in this paper, 
the word ``manifold" always means a manifold with or without boundary, 
as in \cite{Iver}) and $\Psi \subset \Sigma$ is an order filter 
with $\Psi \not \ni \emptyset$. We denote the {\it orientation sheaf} of $X$ 
over $k$ (c.f. \cite[V.\S3]{Iver}) by $\Or_X$. 
Thus we have $or_X[d] \cong \D_X$ in $D^b(\Sh(X))$. Let 
$U := U_\Psi$ be an open subset with the embedding map $j : U \to X$. 
We have $(\Do(Re_\emptyset)_\Psi)^\dagger \cong j_! j^* 
\Do(R e_\emptyset)^\dagger \cong j_! j^* \D_X \cong j_! \D_U \cong (j_! or_U) 
[d]$. Thus $[ \, \Ext^i_R(\, \Do(Re_\emptyset)_\Psi, \, \wc \, ) \, 
]_\emptyset \cong H_\emptyset^{-i+1} (\Do(Re_\emptyset)_\Psi)^\vee
\cong H^{d-i}_c(U, or_U)^\vee$. But we have 
$H^i(U;k) \cong H^{d-i}_c(U, or_U)^\vee$ by the Poincar\'e duality. 
So the equality in Theorem~\ref{cohomology3} actually holds. 
\end{exmp}

For a finite poset $P$, the {\it order complex} 
$\Delta(P)$ is the set of chains of $P$. Recall that a subset $C$ of $P$ 
is a {\it chain} if any two elements of $C$ are comparable. 
Obviously, $\Delta(P)$ is an (abstract) simplicial complex. 
The geometric realization of the order complex $\Delta(\Sigma ')$ of 
$\Sigma ' := \Sigma  \setminus \emptyset$ is homeomorphic to 
the underlying space $X$ of $\Sigma$. 

We say a finite regular 
cell complex $\Sigma $ is {\it Cohen-Macaulay} (resp. {\it Buchsbaum})
if $\Delta(\Sigma ')$ is Cohen-Macaulay (resp. Buchsbaum) over 
$k$ in the sense of \cite[II.\S\S3-4]{St} (resp. \cite[II.\S8]{St}). 
(If $\Sigma$ itself is a simplicial complex, 
we can use $\Sigma$ directly instead of $\Delta(\Sigma')$.) 
These are topological properties of the underlying space $X$. In fact,  
$\Sigma$ is Buchsbaum if and only if $\cH^i(\D_X) = 0$ 
for all $-i \ne d := \dim X$ (see  \cite[Corollary~4.7]{Y6}). 
For example, if $X$ is a manifold, $\Sigma$ is Buchsbaum. 
Similarly, $\Sigma$ is Cohen-Macaulay if and only if it is Buchsbaum and 
$\rH^i(X;k) = 0$ for all $i < d$.

We have  $$H^i(\Do(R e_\emptyset))_\emptyset = 
\Ext^i_R(R e_\emptyset, \wc) \cong 
H_\emptyset^{-i+1}(R e_\emptyset)^\vee \cong \rH^{-i}(X;k)^\vee \quad 
\text{for all $i \in \ZZ$}$$ 
by Corollary~\ref{Serre duality} and Theorem~\ref{cohomology}. 
Recall that $\Do(R e_\emptyset)^\dagger \cong \D_X$. 
So $H^i(\Do(R e_\emptyset)) = 0$ for all 
$i \ne -d$ if and only if $X$ is Cohen-Macaulay over $k$.  
In general,  $H^i(\wc)^\dagger$ can be non-zero for some  $i \ne -d$ 
even if $X$ is Cohen-Macaulay. For example, let $X$ be a closed 2-dimensional 
disc, and $\Sigma$ the regular cell decomposition of $X$ given in 
Example~\ref{2 disc}. Then the ``$\rho_1$-$\rho_2$ component" 
$(\wc)^{\rho_1}_{\rho_2}$ of $\wc$ is of the form 
$0 \to E(\sigma)^{\rho_1}_{\rho_2} \to E(\tau_1)^{\rho_1}_{\rho_2} \oplus 
E(\tau_2)^{\rho_1}_{\rho_2} \to 0$. Thus 
$H^{-1}(\wc)^{\rho_1}_{\rho_2} \ne 0$. But we have the following.

\begin{prop}\label{w = module}
Assume that $\Sigma$ is a meet-semilattice as a poset (e.g., $\Sigma$ is 
a simplicial complex). Then we have the following. 

(1) $H^i(\wc) = 0$ for all $i \ne -d$ if and only if 
$\Sigma$ is Cohen-Macaulay over $k$.

(2) $H^i(\wc)^\dagger = 0$ for all $i \ne -d$ if and only if 
$\Sigma$ is Buchsbaum over $k$.
\end{prop}

\begin{proof}
(1) Since $\wc \cong \Do(R) \cong 
\bigoplus_{\sigma \in \Sigma} \Do(R e_\sigma)$, 
the ``only if part" is clear from the argument before the proposition. 
To prove the ``if part", we assume that $\Sigma$ is Cohen-Macaulay. 
Set $\Omega := H^{-d}(\Do(R e_\emptyset))$. 
Then $\Omega[d] \cong \Do(R e_\emptyset)$ in $D^b(\md)$. 
By Proposition~\ref{uHom and Do}, we have $\Do(R e_\sigma) \cong 
\bR \uHom_R(R e_\sigma, \Omega[d])$. Since $R e_\sigma$ is a projective 
module, we have $\uExt^i_R(R e_\sigma, \Omega) = 0$ for all $i > 0$ by 
Lemma~\ref{RuHom}. Thus $H^i( \Do(R e_\sigma) ) = 0$ for all $i \ne -d$.

(2) Similar to (1). 
\end{proof}

\begin{rem}
By \cite[Proposition~4.10]{Y6}, 
Proposition~\ref{w = module} (1) states  
that if $\Sigma$ is a Cohen-Macaulay simplicial complex, 
the relative simplicial complex $(\Sigma, \del_\Sigma(\sigma))$ 
is Cohen-Macaulay in the sense of \cite[III.\S7]{St} 
for all $\sigma \in \Sigma$.  
Here $\del_\Sigma(\sigma) := \{ \, \tau \in \Sigma \mid 
\tau \not \geq \sigma  \, \}$ is a subcomplex of $\Sigma$. 
\end{rem}

\begin{exmp}
(1) We say a finite regular cell complex $\Sigma$ of dimension $d$ 
is {\it Gorenstein*} over $k$ (see \cite[p.67]{St}), 
if the order complex $\Delta := \Delta(\Sigma ')$ of $\Sigma ' := 
\Sigma  \setminus \emptyset$ is Cohen-Macaulay over $k$ (i.e., 
$\rH_i(\lk_\Delta \sigma ; k) = 0$ for all $\sigma \in \Sigma$ and 
all $i \ne d - \dim \sigma -1$)  
and $\rH_{d - \dim \sigma -1}(\lk_\Delta \sigma ; k) = k$ for all 
$\sigma \in \Delta$. (If $\Sigma$ itself is a simplicial complex, 
we can use $\Sigma$ directly instead of $\Delta$.) 
This is a topological property of the underlying space $X$. 
For example, if $X$ is homeomorphic to a $d$-dimensional sphere, 
then $\Sigma$ is Gorenstein*.  

It is easy to see that $\Do(R e_\emptyset) \cong (R e_\emptyset)[d]$ in 
$D^b(\md)$ if and only if $X$ is Gorenstein*. If $\Sigma$ is a Gorenstein* 
{\it simplicial} complex, then $\wc \cong \Omega[d]$ for some 
$\Omega \in \mod_{R \otimes_k R}$ by Proposition~\ref{w = module}. 
Moreover, we can describe $\Omega$ explicitly. In fact, $\Omega$ has a 
$k$-basis $\{ \, e_\tau^\sigma \mid \sigma, \tau \in \Sigma, 
\sigma \cup \tau \in \Sigma \, \}$ 
and its module structure is defined by 
$$(e_{\sigma', \, \tau'} \otimes 1) \cdot e^\tau_\rho 
= \begin{cases} 
e^\tau_{\sigma'} & 
\text{if $\tau' = \rho$ and $\sigma' \cup \tau \in \Sigma$,} \\
0 & \text{otherwise,}
\end{cases}$$
and
$$
(1 \otimes e_{\sigma', \, \tau'}) \cdot e^\tau_\rho = 
\begin{cases}
e^{\sigma'}_\rho & 
\text{if $\tau'=\tau$ and $\sigma' \cup \rho \in \Sigma$,}\\
0 & \text{otherwise.}
\end{cases}$$
To check this, note that the ``$\tau$-$\rho$ component" 
$(\wc)^\tau_\rho$ of $\wc = \< \, e(\sigma)^\tau_\rho \mid 
\sigma \geq \tau, \rho \, \>$ is isomorphic to 
$\tC_{-n-\bullet}(\lk_\Sigma (\tau \cup \rho))$ 
as a complex of $k$-vector spaces, where 
$\tC_{\bullet}(\lk_\Sigma (\tau \cup \rho))$ is the augmented chain complex 
of $\lk_\Sigma (\tau \cup \rho)$ and $n = \dim (\tau \cup \rho) +1$. 
So the description follows from the Gorenstein* property of $\Sigma$. 
It is easy to see that $\Do(R e_\sigma) \cong 
\< \, e_\tau^\sigma \mid  \tau \in \str_\Sigma \sigma \, \>$. 
So $\bR j_* \D_{U_\sigma} \cong j_* \const_{\overline{U}_\sigma}[d]$, where 
$j: \overline{U}_\sigma \to X$ is the embedding map of the closure 
$\overline{U}_\sigma$ of $U_\sigma$. 

(2) Let $\Sigma$ be a finite simplicial complex of dimension $d$, 
and $V$ the set of its vertices.  
Assume that $\Sigma$ is Gorenstein in the sense of \cite[II.\S5]{St}. 
Then there is a subset $W \subset V$ and a Gorenstein* simplicial complex 
$\Delta \subset 2^{V \setminus W}$ such that $\Sigma = 2^W * \Delta$, 
where ``$*$" stands for simplicial join. 
(The Gorenstein property depends on the particular simplicial decomposition 
of $X$.) 
Since a Gorenstein simplicial complex is Cohen-Macaulay, there is $\Omega \in 
\mod_{R \otimes_k R}$ such that $\wc \cong \Omega[d]$. By the argument 
similar to (1), $\Omega$ has a $k$-basis $\{ \, e_\tau^\sigma \mid \sigma 
\cup \tau \in \Sigma, \, \sigma \cup \tau \supset W \, \}$ and its left 
$R \otimes_k R$-module structure is given by the similar way to (1). 

Assume that $\Sigma$ is the $d$-simplex $2^V$. 
Then $\Sigma$ is Gorenstein and $\Omega$ has a $k$-basis 
$\{ \, e_\tau^\sigma \mid \sigma \cup \tau = V \, \}$. Moreover, 
we have a ring isomorphism given by $\varphi: R \ni e_{\sigma, \tau} \mapsto 
e_{\tau^\cmpl, \sigma^\cmpl} \in R^\op$, where $R^\op$ is the opposite ring 
of $R$, and $\sigma^\cmpl := V \setminus \sigma$.  
Thus $R$ has a left $(R \otimes_k R)$-module structure given by 
$(x \otimes y) \cdot r = x \cdot r \cdot \varphi(y)$. 
Then a map given by $R \ni e_{\sigma, \tau} \mapsto e_\sigma^{\tau^\cmpl} 
\in \Omega$ is an isomorphism of $(R \otimes R)$-modules. 
So $R$ is an Auslander regular ring in this case. See \cite[Remark~3.3]{Y5}. 

(3) Assume that $\Sigma$ is a simplicial complex and  
$X$ is a $d$-dimensional manifold which is orientable 
(i.e., $or_X \cong \const_X$) and connected. 
Then $H^i(\wc)^\dagger = 0$ for all $i \ne -d$. It is easy to see that 
$\Omega := H^{-d}(\wc) \in \mod_{R \otimes_k R}$ has a $k$-basis 
$\{ \, e_\tau^\sigma \mid \sigma \cup \tau \in \Sigma \, \}$ and 
the module structure is give by the same way as (1).  
\end{exmp}

\section{The M\"obius function of the poset $\hSig$}
The {\it M\"obius function} of a finite poset $P$ is a function 
$\mu : \{ \, (x,y) \mid \text{$x \leq y$ in $P$} \, \} \to \ZZ$ defined 
by the following way
$$\mu(x,x) = 1 \quad \text{for all $x \in P$} \quad \ \text{and} \quad \ 
\mu(x,y) = - \sum_{x\leq z < y} \mu(x,z) \quad 
\text{for all $x < y$ in  $P$.}$$
See \cite[Chapter 3]{StEn} for a general theory of this function. 

For a finite regular cell complex $\Sigma$, let $\hSig$ be the poset 
obtained from $\Sigma$ adjoining the greatest element $\h1$ 
(even if $\Sigma$ already possess the greatest element, we add 
the new one). 
Then the M\"obius function $\mu$ of $\hSig$ has a topological 
meaning. For example, we have $\mu(\emptyset, \h1) = \rXi(X)$, 
where $\rXi(X)$ is the reduced Euler characteristic 
$\sum (-1)^i \dim_K \rH^i(X;k)$ of $X$. 
When the underlying space $X$ is a manifold, the M\"obius function 
of $\hSig$ is completely determined in \cite[Proposition~3.8.9]{StEn}. 
Here we study the general case. 

For $\sigma \in \Sigma$ with $\dim \sigma > 0$, 
$\{ \, \sigma' \in \Sigma \mid \sigma' < \sigma \, \}$ 
is a regular cell decomposition of $\bar{\sigma} -\sigma$ which is 
homeomorphic to a sphere of dimension $\dim \sigma -1$. 
Hence we have $\mu(\tau, \sigma) = (-1)^{l(\tau, \sigma)}$ 
for $\tau \in \Sigma$ with $\tau \leq \sigma$ 
by \cite[Proposition~3.8.9]{StEn}, 
where $l(\tau, \sigma) := \dim \sigma - \dim \tau$. 
So it remains to determine $\mu(\sigma,\h1)$ for $\sigma \ne \emptyset$. 

\begin{prop}\label{Mobius}
For a cell $\emptyset \ne \sigma \in \Sigma$ 
with $j := \dim \sigma$, we have 
$$\mu(\sigma,\h1) 
= \sum_{i \geq j} (-1)^{i-j+1} \dim_K H_c^i(U_\sigma ; k).$$
Here $H_c^i(U_\sigma ; k)$ is the cohomology with compact support
of the open set $U_\sigma = \bigcup_{\rho \geq \sigma} \rho$ of $X$. 
\end{prop}

\begin{proof}
The assertion follows from the next computation. 
\begin{eqnarray*}
\mu(\sigma,\h1) &=& - \sum_{\rho \in \Sigma, \, 
\rho \geq \sigma} \mu(\sigma, \rho)\\ 
&=& \sum_{i \geq j} (-1)^{i-j+1} \cdot \# \{ \, \rho \in \Sigma \mid 
\rho \geq \sigma, \,  \dim \rho =i \, \} \\ 
&=& \sum_{i \geq j} (-1)^{i-j+1} \dim_K \cH^{-i}(\D_X)(U_\sigma) \\
&=& \sum_{i \geq j} (-1)^{i-j+1} \dim_K H_c^i(U_\sigma ; k). 
\end{eqnarray*}
The second equality follows the fact that 
$\mu(\sigma, \rho) = (-1)^{l(\sigma, \rho)}$. 
The third equality follows from $\D_X \cong 
\Do(R e_\emptyset)^\dagger$ and the description of $\Do(R e_\emptyset)$. 
Recall also that $M^\dagger(U_\sigma) \cong M_\sigma$. 
And the last equality follows from the Verdier duality. 
\end{proof}

Assume that $X$ is a manifold of dimension $d$.  
If $\sigma \ne \emptyset$ is contained in the boundary of $X$, 
then $U_\sigma$ is homeomorphic to $(\RR^{d-1} \times \RR_{\geq 0})$ and 
$H_c^i(U_\sigma ; k)=0$ for all $i$. Thus $\mu(\sigma, \h1) = 0$ in this 
case. If $\sigma$ is not contained in the boundary of $X$, 
then $U_\sigma$ is homeomorphic to $\RR^d$ and 
$H_c^i(U_\sigma ; k)=0$ for all $i \ne d$ and $H_c^d(U_\sigma ; k)=k$. 
Hence we have $\mu(\sigma, \h1) = (-1)^{d -\dim \sigma +1}$. 
So Proposition~\ref{Mobius} recovers \cite[Proposition~3.8.9]{StEn}.

\section{Relation to Koszul duality}
Let $A = \bigoplus_{i \geq 0} A_i$ be an $\NN$-graded associative $k$-algebra 
such that $\dim_k A_i < \infty$ for all $i$ and $A_0 \cong k^n$ for some 
$n \in \NN$ as an algebra. Then $\rad := \bigoplus_{i > 0} A_i$ is the graded 
Jacobson radical.  We say $A$ is {\it Koszul}, 
if a left $A$-module $A/\rad$ admits a graded projective resolution 
$$\cdots \to P^{-2} \to P^{-1} \to P^0 \to A/\rad \to 0$$ such that $P^{-i}$ 
is generated by its degree $i$ component as an $A$-module 
(i.e., $P^{-i} = A P^{-i}_i$). 
If $A$ is Koszul, it is a quadratic ring, and its {\it quadratic dual 
ring} $A^!$ (see \cite[Definition~2.8.1]{BGS}) 
is Koszul again, and isomorphic to the opposite ring of the Yoneda algebra 
$\Ext^\bullet_A(A/\rad, A/\rad)$. 

\medskip

Note that the incidence algebra $R$ of $\Sigma $  is a graded ring 
with $\deg(e_{\sigma,\sigma'}) = \dim \sigma - \dim \sigma'$. 
So we can discuss the Koszul property of $R$.

\begin{prop}[c.f. {\cite[Lemma~4.5]{Y5}}]\label{self dual}
The incidence algebra $R$ of a finite regular cell complex $\Sigma$ 
is always Koszul.  
And the quadratic dual ring  $R^!$ is isomorphic to $R^\op$. 
\end{prop}

When $\Sigma$ is a simplicial complex,  
the above result was proved by Polishchuk \cite{Pls} in much wider 
context(but, $\emptyset \not \in \Sigma$ in his convention). 
More precisely, he put the new partial order on the set 
$\Sigma \setminus \emptyset$ associated with a perversity function $p$, 
and construct two rings from this new poset. 
Then he proved that these two rings are Koszul and quadratic dual rings of 
each other. Our $R$ and $R^\op$ correspond to the case 
when $p$ is a bottom (or top) perversity. 
In the middle perversity case, $\Sigma$ has to be a {\it simplicial} complex 
to make their rings Koszul. 

\begin{proof}
By \cite{P,W}, $R$ is Koszul if and only if the order complex $\Delta(I)$ 
is Cohen-Macaulay over $k$ for any open interval $I$ of $\Sigma $. 
Set $\Sigma ' := \Sigma  \setminus \emptyset$. 
Note that $\Delta(I) = \lk_{\Delta(\Sigma')}F$ for some $F 
\in \Delta(\Sigma ')$ containing a maximal cell $\sigma \in \Sigma$. 
Set $\Delta := \str_{\Delta(\Sigma')}\sigma$. Then $\Delta(I) = \lk_\Delta F$. 
Since  the underlying space of $\Delta$ is the closed disc 
$\bar{\sigma}$, $\Delta$ is Cohen-Macaulay. 
Hence $\lk_\Delta F$ is also. So $R$ is Koszul. 

Let $T := T_{R_0} R_1 = R_0 \oplus R_1 \oplus 
(R_1 \otimes_{R_0} R_1) \oplus \cdots 
= \bigoplus_{i \geq 0} 
R_1^{\otimes i}$ be the tensor ring of $R_1 = 
\< \, e_{\sigma, \, \tau} \mid \sigma,\tau \in \Sigma , \, 
\sigma > \tau, \, \dim \sigma  = \dim \tau +1  \, \>$ over $R_0$. 
Then $R \cong T/I$, where 
$$I = (\,  e_{\sigma, \, \rho_1} \otimes e_{\rho_1, \, \tau}  - 
e_{\sigma, \, \rho_2} \otimes e_{\rho_2, \, \tau} \mid \sigma, 
\tau, \rho_1, \rho_2 \in \Sigma , \, \sigma > \rho_i > \tau, \, 
\dim \sigma  = \dim \tau +2 \,)$$ 
is a two sided ideal.  Let $R_1^* := \Hom_{R_0}(R_1, R_0)$ 
be the dual of the {\it left} $R_0$-module $R_1$. 
Then  $R_1^*$ has a right $R_0$-module 
structure such that $(fa)(v) = (f(v))a$, and a left $R_0$-module 
structure such that $(af)(v) =  f(va)$, where $a \in R_0$, 
$f \in R_1^*$, $v \in R_1$. As a left (or right) 
$R_0$-module, $R_1^*$ is generated by  
$\{ \, e_{\tau,\, \sigma}^* \mid  
\sigma > \tau, \, \dim \sigma  = \dim \tau +1 \, \}$, where 
$e_{\tau, \, \sigma}^*(e_{\sigma', \, \tau'}) = 
\delta_{\sigma, \sigma'} \cdot \delta_{\tau, \tau'} \cdot e_\sigma.$ 

Let $T^* = T_{R_0} R_1^*$ be the tensor ring of 
$R_1^*$. Note that $e^*_{\tau, \, \sigma} \otimes 
e^*_{\tau', \, \sigma'} \in R_1^* \otimes_{R_0} 
R_1^*$ is non-zero if and only if $\sigma = \tau'$. 
We have that $(R_1^* \otimes_{R_0} R_1^*)$ 
is isomorphic to $(R_1 \otimes_{R_0} R_1)^* = 
\Hom_{R_0}(R_1 \otimes_{R_0} R_1, R_0)$ 
via $(f \otimes g)(v \otimes w) = g(vf(w))$, where 
$f, g \in R_1^*$ and $v, w \in R_1$. 
In particular, 
$( e^*_{\tau, \, \rho} \otimes e^*_{\rho, \, \sigma}) 
( e_{\sigma, \, \rho} \otimes e_{\rho, \, \tau} ) = e_\sigma.$
Recall that if $\sigma, \tau \in \Sigma $, $\sigma > \tau$ 
and $\dim \sigma  = \dim \tau +2$, 
then there are exactly two cells $\rho_1, \rho_2 \in \Sigma $ 
between $\sigma$ and $\tau$. 
So easy computation shows that the quadratic dual ideal 
$$I^\bot = (\, f \in R_1^* \otimes R_1^* 
\mid \text{$f(v) = 0$ 
for all $v \in I_2 \subset  R_1 \otimes R_1 = T_2$} \,) 
\subset T^*$$ of $I$ is equal to  
$$(\, e^*_{\tau, \, \rho_1} \otimes e^*_{\rho_1, \, \sigma}  + 
e^*_{\tau, \, \rho_2} \otimes e^*_{\rho_2, \, \sigma} \mid 
\sigma, \tau, \rho_1, \rho_2 \in \Sigma , \, \sigma > \rho_i > 
\tau, \, \dim \sigma  = \dim \tau +2 \, ).$$ 
The  $k$-algebra homomorphism $R \to R^! = T^*/I^\bot$  
defined by the identity map on 
$R_0 = T_0 = (T^*)_0 = (R^!)_0$ and $R_1 \ni e_{\sigma, \, \tau} 
\mapsto \varepsilon(\sigma,\tau) \cdot e^*_{\tau, \, \sigma} 
\in R^!_1$ is a graded isomorphism. Here $\varepsilon$ 
is a incidence function of $\Sigma $. 
\end{proof}

Since $R^! \cong R^\op$,  
$\Hom_k(-,k)$ gives duality functors $\Dk: \md \to \mdop$ and 
$\Dko: \mdop \to \md$. These functors are exact, and they can be extended 
to the duality functors between $D^b(\md)$ and $D^b(\mdop)$. 

Note that $R^!$ is a graded ring with 
$\deg e^*_{\tau, \, \sigma} = \dim \sigma - \dim \tau$. 
Let $\mdZ$ (resp. $\mdZop$) be the category of finitely generated graded left 
$R$-modules (resp. $R^!$-modules). 
Note that we can regard the functor $\Do$ (resp. $\Dk$ and $\Dko$)
as the functor from $D^b(\mdZ)$ to itself (resp. $D^b(\mdZ) \to D^b(\mdZop)$ 
and $D^b(\mdZop) \to D^b(\mdZ)$). 

For each $i \in \ZZ$, let $\mdZ(i)$ be the full 
subcategory of $\mdZ$ consisting of $M \in \mdZ$ with 
$\deg M_\sigma = \dim \sigma -i$. 
For any $M \in \mdZ$, there are modules  $M^{(i)} \in \mdZ(i)$ such that 
$M \cong \bigoplus_{i \in \ZZ}M^{(i)}$. 
The forgetful functor gives an equivalence 
$\mdZ(i) \cong \md$ for all $i \in \ZZ$, and $D^b(\mdZ(i))$ 
is a full subcategory of $D^b(\mdZ)$.  
Similarly, let $\mdZop(i)$ be the full 
subcategory of $\mdZop$ consisting of $M \in \mdZ$ with $\deg M_\sigma 
= -\dim \sigma-i$. The above mentioned facts on 
$\mdZ(i)$ also hold for $\mdZop(i)$.

Let $DF: \md \to \mdop$ and $DG: \mdop \to \md$ be the functors defined 
in \cite[Theorem~2.12.1]{BGS}. 
Since $R$ and $R^!$ are artinian, 
$DF$ and $DG$ give an equivalence $D^b(\mdZ) \cong D^b(\mdZop)$
by the Koszul duality (\cite[Theorem~2.12.6]{BGS}). 

When $\Sigma$ is a simplicial complex 
the next result was given by Vybornov \cite{Vyb} 
(under the convention that  $\emptyset \not \in \Sigma$). Independently, 
the author also proved a similar result (\cite[Theorem~4.7]{Y5}).

\begin{thm}[c.f. Vybornov, {\cite[Corollary~4.3.5]{Vyb}}]\label{DF, DG}
Under the above notation, if $\M \in D^b(\mdZ(0))$,  then 
we have $DF(\M) \in D^b(\mdZop(0))$. Similarly, 
if $\N \in D^b(\mdZop(0))$,  then $DG(\N) \in D^b(\mdZ(0))$. Under 
the equivalence $\mdZ(0) \cong \md$ and $\mdZop(0) \cong \mdop$, we have 
$DF \cong \Dk \circ \Do$ and $DG \cong \Do \circ \Dko$. 
\end{thm}

\begin{proof}
 Recall that $(R^!)_0 = R_0$. Let $N \in \mdop$.
For the functor $DG$, we need the left 
$R$-module structure on $\Hom_{R_0}(R, N_\sigma)$ 
given by $(xf)(y) :=f(yx)$. The $R$-morphism given by 
$\Hom_{R_0}(R, N_\sigma) \ni f \longmapsto 
\sum_{\tau \leq \sigma} e(\sigma)_\tau \otimes_k f(e_{\sigma, \, \tau}) 
\in E(\sigma) \otimes_k N_\sigma$ 
gives an isomorphism $\Hom_{R_0}(R, N_\sigma) \cong
E(\sigma) \otimes_k N_\sigma$. Under this isomorphism, for cells 
$\tau < \sigma$, the morphism 
$\Hom_{R_0}(R, N_\sigma) \to 
\Hom_{R_0}(R, N_\tau)$ given by 
$f \mapsto [x \mapsto e^*_{\tau, \, \sigma} \, f(e_{\sigma, \, \tau} \,x)]$ 
corresponds to the morphism 
$E(\sigma) \otimes_k N_\sigma \to E(\tau) \otimes_k N_\tau$ given by 
$e(\sigma)_\rho \otimes y \mapsto e(\tau)_\rho \otimes 
e^*_{\tau,\,\sigma} \, y$. 
(Here  $e(\tau)_\rho = 0$ if $\tau \not \geq \rho$.) 

Let $N \in \mdZop$. By the explicit description of $\Do$ given in \S3, we have 
$$(\Do \circ \Dko)^i(N) = 
\bigoplus_{\substack{\sigma \in \Sigma  \\ \dim \sigma = -i}} 
E(\sigma)  \otimes_k N_\sigma\\
= \bigoplus_{\substack{\sigma \in \Sigma  \\ \dim \sigma = -i}}
\Hom_{R_0}(R, N_\sigma)$$
and the differential map defined by 
$$E(\sigma)  \otimes_k N_\sigma  \ni e(\sigma)_\rho \otimes y \mapsto 
\sum_{\substack{\tau \in \Sigma  \\ \dim \tau = -i-1}}
\varepsilon(\sigma,\tau) \, ( \, e(\tau)_\rho
\otimes e^*_{\tau, \, \sigma} \, y) 
\in (\Do \circ \Dko)^{i+1}(N).$$
So, if we forget the grading of modules, 
we have $DG(N) \cong (\Do \circ \Dko) (N)$. 
Similarly, we can check an isomorphism 
$DG(\N) \cong (\Do \circ \Dko) (\N)$ 
for a complex $\N \in D^b(\mdZop)$. 

Assume that $N \in \mdZop(0)$. Then the degree of 
$e(\sigma)_\tau \otimes y \in E(\sigma) \otimes_k N_\sigma \subset DG(N)$ 
is $(\dim \tau - \dim \sigma)+\dim \sigma = \dim \tau$ 
(see the proof of \cite[Theorem~2.12.1]{BGS} for the grading 
of $DG(N)$). Thus we have $DG(N) \in \mdZ(0)$.  

We can prove the statement on $DF$ by a similar (easier) way. 
\end{proof}

The results corresponding to Proposition~\ref{self dual} 
and Theorem~\ref{DF, DG} also hold for the incidence algebra of the poset 
$\Sigma \setminus \emptyset$. In other words, 
Vybornov~\cite[Corollary~4.3.5]{Vyb} and the ``top perversity case" of 
Polishchuk~\cite{Pls} can be generalized directly 
into regular cell complexes.

\section*{Acknowledgments}
The author is grateful to Professor Maxim Vybornov. 
He informed me his works (including \cite{Vyb}) and related papers 
(including \cite{Pls}), and sent me Shepard's thesis \cite{Sh}. 
The author has learned a lot of things from theses papers.

\end{document}